\documentclass[preprint]{elsarticle}
\usepackage{amsmath}
\usepackage{amssymb}
\usepackage{amsfonts}
\usepackage{graphicx}
\usepackage{lineno,hyperref}
\usepackage{subfigure}
\newtheorem{thm}{Theorem}
\newtheorem{lem}[thm]{Lemma}

\newtheorem{corollary}[thm]{Corollary}

\newdefinition{rmk}{Remark}
\newdefinition{definition}{Definition}
\newproof{pf}{Proof}

\modulolinenumbers[5]

\journal{Journal of Geometry and Physics}

\begin{document}

\begin{frontmatter}

\title{On equivalence of the second order linear differential operators,
acting in vector bundles}

\author[mymainaddress]{Valentin Lychagin\corref{mycorrespondingauthor}}
\cortext[mycorrespondingauthor]{Corresponding author}
\ead{valentin.lychagin@uit.no}

\address[mymainaddress]{V.A. Trapeznikov Institute of Control Sciences, Russian Academy of Sciences, 65 Profsoyuznaya Str., 117997 Moscow, Russia}

\begin{abstract}
The equivalence problem for linear differential operators of the second order, acting in vector bundles, is discussed. The field of rational
invariants of symbols is described and connections, naturally associated
with differential operators, are found. These geometrical structures are
used to solve the problems of local as well as global equivalence of
differential operators.
\end{abstract}


\end{frontmatter}

\section{Introduction}

This paper is the last in the series of papers devoted to invariants,
geometrical structures and equivalence of linear differential operators,
acting in vector bundles (\cite{Ly2},\cite{LyB},\cite{LYk}). \ Here we study
linear differential operators of the second order, acting in sections of a
vector bundle $\pi :E\left( \pi \right) \rightarrow M,$ where $\dim M=n\geq
2,\dim \pi =m\geq 2.$

Our approach is based on the following main steps. First of all we study
geometrical structures defined by symbols of differential operators and use
them to define invariants and orbits of the symbols with respect to groups,
induced by automorphisms of bundles. The second step is extremely important
for us and consist of finding connections, naturally associated with
differential operators. Thus in the given case, in general, the symbol of
the differential operators defines a pseudo Riemannian structure on the base
manifold and the condition of the triviality of the subsymbol leads us to
uniquely defined connection in the bundle. This connection is naturally
related to the differential operator and gives us the quantization, or
splitting of differential operator into the sum of the symbol and an
endomorphism of the bundle. After this, the splitting and the Artin-Procesi
theorem (\cite{Art},\cite{Pr}) allow us to find invariants of operators. The
last step of the approach (similar to the one used in \cite{Ly2},\cite{LyB},\cite%
{LYk}) consist of in using the invariants to construct the natural
coordinates and local models of differential operators. This gives us local
and global description of differential operators under some natural
conditions of regularity on the operator.

\section{Differential operators}

\subsection{Preliminaries}

The notations we use in this paper are similar to the ones that were
used in the papers \cite{LyB} and \cite{LYk} .

Let $M$ be an $n$-dimensional manifold and let $\pi :E\left( \pi \right)
\rightarrow M$ be a vector bundle.

We denote by $\tau :TM\rightarrow M$ and $\tau ^{\ast }:T^{\ast
}M\rightarrow M$ \ the tangent and respectively cotangent bundles over the
manifold $M,$ and by $\mathbf{1:}\mathbb{R}\times M\rightarrow M\ $we denote
the trivial linear bundle $.$

The symmetric and exterior powers of a vector bundle $\pi :E\left( \pi
\right) \rightarrow M$ \ will be denoted by $\mathbf{S}^{k}\left( \pi
\right) $ and $\mathbf{\Lambda }^{k}\left( \pi \right) .$

The module of smooth sections of bundle $\pi $ we denote by $C^{\infty
}\left( \pi \right) ,$ and for the cases tangent, cotangent and the trivial
bundles we will use the following notations: $\Sigma _{k}\left( M\right)
=C^{\infty }\left( \mathbf{S}^{k}\left( \tau \right) \right) -$ is the module
of symmetric $k$-vectors and $\Sigma ^{k}\left( M\right) =C^{\infty }\left(
\mathbf{S}^{k}\left( \tau ^{\ast }\right) \right) -$ is the module of symmetric
$k$-forms, $\Omega _{k}\left( M\right) =C^{\infty }\left( \mathbf{\Lambda }%
^{k}\left( \tau \right) \right) -$ is the module of skew-symmetric $k$-vectors
and
 $\Omega ^{k}\left( M\right) =C^{\infty }\left( \mathbf{\Lambda }%
^{k}\left( \tau ^{\ast }\right) \right) $ -is the module of exterior $k$-forms,
$C^{\infty }\left( \mathbf{1}\right) =C^{\infty }\left( M\right) .$

Let $\mathbf{Diff}_{k}\mathbf{(\pi ,\pi )}$ be the module of linear
differential operators of order $k$ acting in the sections of the bundle $%
\pi .$ These modules are connected by the exact sequences
\begin{equation}
\mathbf{0\rightarrow Diff}_{k-1}\left( \pi ,\pi \right) \mathbf{\rightarrow
Diff}_{k}\left( \pi ,\pi \right) \overset{\mathrm{smbl}_{k}}{\rightarrow }%
\mathrm{End}\left( \pi \right) \otimes \Sigma _{k}\left( M\right)
\rightarrow \mathbf{0,}  \label{exact1}
\end{equation}%
where $\mathrm{End}\left( \pi \right) =\mathbf{Diff}_{0}\left( \pi ,\pi
\right) $ is the module of endomorphisms of $C^{\infty }\left( \pi \right) ,$
the mapping $\mathrm{smbl}$ assigns to differential operator $\Delta $ its
symbol $\mathrm{smbl}_{k}\left( \Delta \right) $ and the tensor product
 is taking over $ C^{\infty }\left( M\right) .$

If we consider $\mathrm{smbl}_{k}\left( \Delta \right) $ as a linear
operator $\mathrm{smbl}_{k}\left( \Delta \right) :\Sigma ^{k}\left(
M\right) \rightarrow \mathrm{End}\left( \pi \right) ,$ then
\begin{equation*}
\mathrm{smbl}_{k}\left( \Delta \right) \left( df^{k}\right) =\frac{%
1}{k!}\delta _{f}^{k}\left( \Delta \right) \in \mathrm{End}\left( \pi
\right) ,
\end{equation*}%
where $f\in $ $C^{\infty }\left( M\right) ,\ \delta _{f}:\mathbf{Diff}%
_{i}\left( \pi ,\pi \right) \mathbf{\rightarrow Diff}_{i-1}\left( \pi ,\pi
\right) ,$ is the commutator mapping $\delta _{f}\left( \square \right)
=f\circ \square -\square \circ f,$ and $df^{k}=df\cdot \cdots \cdot df\in
\Sigma ^{k}\left( M\right) $ is the $k$-th symmetric power of the
differential $1$-form $df\in \Sigma ^{1}\left( M\right) .$

\subsection{Pseudogroup actions}

As before (\cite{LyB}), we consider two pseudogroups: $\mathcal{G}\left(
M\right) -$ the pseudogroup of local diffeomorphisms of manifold $M,$ and $%
\mathbf{Aut}(\pi )-$ the pseudogroups of local automorphisms of vector
bundle $\pi $ over $M.$

The following sequence of pseudogroup morphisms is exact
\begin{equation}
1\rightarrow \mathrm{GL}\left( \pi \right) \rightarrow \mathbf{Aut}(\pi
)\rightarrow \mathcal{G}\left( M\right) \rightarrow 1,
\label{group exact seq}
\end{equation}%
where $\mathrm{GL}\left( \pi \right) \subset \mathbf{Aut}(\pi )$ is the
pseudogroup of automorphisms that are identity on $M.$

We will consider the natural actions of these pseudogroups on sections of the
bundles and on operators. Namely, let $\widetilde{\phi }$ be a local
automorphism, $\widetilde{\phi }\in \mathbf{Aut}(\pi ),$ covering a local
diffeomorphism $\phi \in \mathcal{G}\left( M\right) .$

We define action of $\widetilde{\phi }$ on sections $s\in C^{\infty }\left(
\pi \right) $ as
\begin{equation*}
\widetilde{\phi }_{\ast }:s\longmapsto \widetilde{\phi }\circ s\circ \phi
^{-1},
\end{equation*}%
and
\begin{equation*}
\widetilde{\phi }_{\ast }:\Delta \longmapsto \widetilde{\phi }_{\ast }\circ
\Delta \circ \widetilde{\phi _{\ast }}^{-1},
\end{equation*}%
for differential operators.

\section{Quantizations, symbols and subsymbols}

Let $\Sigma ^{\cdot }=\oplus _{k\geq 0}\Sigma ^{k}\left( M\right) $ be the
graded algebra of symmetric differential forms and let $\Sigma ^{\cdot
}\left( \pi \right) =C^{\infty }\left( \pi \right) \otimes \Sigma ^{\cdot }$
be the graded $\Sigma ^\mathbf{\cdot }$-module of symmetric differential forms with
values in the bundle $\pi .$

Let $\nabla $ be a connection in the vector bundle $\pi $ and let $\nabla
_{c}$ be a connection in the cotangent bundle $\tau ^{\ast }.$

Then their covariant differentials
\begin{equation*}
d_{\nabla _{c}}:\Omega ^{1}\left( M\right) \rightarrow \Omega ^{1}\left(
M\right) \otimes \Omega ^{1}\left( M\right) ,
\end{equation*}%
and
\begin{equation*}
d_{\nabla }:C^{\infty }\left( \pi \right) \rightarrow C^{\infty }\left( \pi
\right) \otimes \Omega ^{1}\left( M\right)
\end{equation*}%
define derivation%
\begin{equation*}
d_{\nabla _{c}}:\Sigma ^{\cdot }\rightarrow \Sigma ^{\cdot +1},
\end{equation*}%
of degree one in graded algebra $\Sigma ^{\cdot }$ and derivation
\begin{equation*}
d_{\nabla }:\Sigma ^{\cdot }\left( \pi \right) \rightarrow \Sigma ^{\cdot
+1}\left( \pi \right)
\end{equation*}%
of degree one in graded $\Sigma ^{\cdot }$-module $\Sigma ^{\cdot }\left(
\pi \right) $ respectively.

Indeed, derivations are completely defined by their actions on generators.

In our case, we get%
\begin{eqnarray*}
d_{\nabla _{c}} &=&d:C^{\infty }\left( M\right) \rightarrow \Omega
^{1}\left( M\right) =\Sigma ^{1}, \\
d_{\nabla _{c}} &:&\Omega ^{1}\left( M\right) =\Sigma ^{1}\overset{d_{\nabla
}}{\rightarrow }\Omega ^{1}\left( M\right) \otimes \Omega ^{1}\left(
M\right) \overset{\mathrm{Sym}}{\rightarrow }\Sigma ^{2},
\end{eqnarray*}%
where $\mathrm{Sym}$ denotes the symmetrization, and define $d_{\nabla }$ as a derivation over $d_{\nabla _{c}}$ such that
\begin{equation*}
d_{\nabla }:C^{\infty }\left( \pi \right) \rightarrow C^{\infty }\left( \pi
\right) \otimes \Sigma ^{1}.
\end{equation*}%
is the covariant differential.

Let now $\sigma \in \mathrm{End}\left( \pi \right) \otimes \Sigma _{k}$
be a symbol. We define the differential operator $\widehat{\sigma }\in \mathbf{%
Diff}_{k}\left( \pi ,\pi \right) $ as follows:%
\begin{equation}
\widehat{\sigma }\left( s\right) \overset{\text{def}}{=}\frac{1}{k!}%
\left\langle \sigma ,\left( d_{\nabla }\right) ^{k}\left( s\right)
\right\rangle ,
\end{equation}

where $s\in C^{\infty }\left( \pi \right) ,\left( d_{\nabla }\right)
^{k}\left( s\right) \in C^{\infty }\left( \pi \right) \otimes \Sigma ^{k},$%
and $\left\langle \cdot ,\cdot \right\rangle $ is the natural pairing
\begin{equation*}
\mathrm{End}\left( \pi \right) \otimes \Sigma _{k}\otimes C^{\infty }\left(
\pi \right) \otimes \Sigma ^{k}\rightarrow C^{\infty }\left( \pi \right) .
\end{equation*}%
Remark that the value of the symbol of the derivation $d_{\nabla }$ on a
covector $\theta $ equals the symmetric product by $\theta $ in the
module $\Sigma ^{\cdot }\left( \pi \right) $ and because the symbol of a
composition of operators equals the composition of symbols we get that the
symbol of the operator $\widehat{\sigma }$ equals $\sigma .$

We call this operator $\widehat{\sigma }$  \textit{the quantization of symbol }$%
\sigma $ and write $\ \widehat{\sigma }=Q_{\nabla }\left( \sigma \right) .$

By the construction morphism $Q_{\nabla }:\mathrm{End}\left( \pi \right)
\otimes \Sigma _{k}\rightarrow \mathbf{Diff}_{k}\left( \pi ,\pi \right) $
splits sequence (\ref{exact1}).

Let now $A\in \mathbf{Diff}_{k}\left( \pi ,\pi \right) $ be a differential
operator and $\sigma _{k}\left( A\right) \in \mathrm{End}\left( \pi \right)
\otimes \Sigma _{k}$ be its symbol. Then operator
\begin{equation*}
A-Q_{\nabla }\left( \sigma _{k}\left( A\right) \right)
\end{equation*}%
has order $\left( k-1\right) ,$ and let $\sigma _{k-1}\left( A,\nabla
\right) \in \mathrm{End}\left( \pi \right) \otimes \Sigma _{k-1}$ be its
symbol.

Then operator $A-Q_{\nabla }\left( \sigma _{k}\left( A\right) \right)
-Q_{\nabla }\left( \sigma _{k-1}\left( A,\nabla \right) \right) $ has order $%
\left( k-2\right) .$ Repeating this process we get subsymbols $\sigma
_{i}\left( A,\nabla \right) \in \mathrm{End}\left( \pi \right) \otimes
\Sigma _{i},$ $0\leq i\leq k-1,$ such that
\begin{equation*}
A=Q_{\nabla }\left( \sigma \left( A,\nabla \right) \right) ,
\end{equation*}%
where
\begin{equation*}
\sigma \left( A,\nabla \right) =\oplus _{0\leq i\leq k}\sigma _{i}\left(
A,\nabla \right)
\end{equation*}%
is a \textit{total symbol }of the operator, and $Q_{\nabla }\left( \sigma
\left( A,\nabla \right) \right) =\sum_{i}Q_{\nabla }\left( \sigma _{i}\left(
A,\nabla \right) \right) .$

\begin{thm}
Let $\ \widetilde{\nabla }$ be a connection in the bundle $\pi $ and let $\
\widetilde{\nabla }=\nabla +\alpha \ \ $be another connection, where $\alpha
\in \mathrm{End}\left( \pi \right) \otimes \Sigma ^{1}\left( M\right) ,$ $%
\widetilde{\nabla }_{X}=\nabla _{X}+\alpha \left( X\right) .$ \newline
Then
\begin{equation*}
\sigma _{k-1}\left( A,\nabla \right) =\sigma _{k-1}\left( A,\widetilde{%
\nabla }\right) +\left\langle \sigma _{k},\alpha \right\rangle ,
\end{equation*}%
where $\left\langle \sigma _{k},\alpha \right\rangle \in \mathrm{End}\left(
\pi \right) \otimes \Sigma _{k-1}$ is the natural pairing $\mathrm{End}%
\left( \pi \right) \otimes \Sigma _{k}\times \mathrm{End}\left( \pi \right)
\otimes \Sigma ^{1}\left( M\right) \rightarrow \mathrm{End}\left( \pi
\right) \otimes \Sigma _{k-1}.$
\end{thm}

\begin{pf}
Remark that $\delta _{f}\left( d_{\nabla }\right) =df$ \ and
\begin{equation*}
\delta _{f}^{i}\left( d_{\nabla }^{k}\right) =k\left( k-1\right) \cdots
\left( k-i+1\right) \ df^{i}\cdot d_{\nabla }^{k-i}.
\end{equation*}%
\newline
Therefore,
\begin{equation*}
\delta _{f}^{k-1}\left( d_{\nabla }^{k}\right) =k!df^{k-1}\cdot d_{\nabla },
\end{equation*}%
and
\begin{equation*}
\delta _{f}^{k-1}\left( Q_{\nabla }\left( \sigma _{k}\right) \right) \left(
s\right) =\left\langle \sigma _{k},df^{k-1}\cdot d_{\nabla }\left( s\right)
\right\rangle .
\end{equation*}%
\newline
Applying this formula to connection $\nabla +\alpha $ we get the statement
of the theorem:%
\begin{eqnarray*}
&&df^{k-1}\rfloor \left( \sigma _{k-1}\left( A,\nabla \right) -\sigma
_{k-1}\left( A,\widetilde{\nabla }\right) \right) =\delta _{f}^{k-1}\left(
Q_{\widetilde{\nabla }}\left( \sigma _{k}\right) -Q_{\nabla }\left( \sigma
_{k}\right) \right) = \\
&&\left\langle \sigma _{k},df^{k-1}(d_{\nabla +\alpha }-d_{\nabla
})\right\rangle =\left\langle \sigma _{k},df^{k-1}\cdot \alpha \right\rangle
=df^{k-1}\rfloor \left\langle \sigma _{k},\alpha \right\rangle ,
\end{eqnarray*}
\end{pf}
where $ \rfloor $ stands for the inner product.
We apply now this result to operators of the order two.

At first, we say that the symbol $\sigma _{2}\in \mathrm{End}\left( \pi
\right) \otimes \Sigma _{2}$ of an operator $A\in \mathbf{Diff}_{2}\left(
\pi ,\pi \right) $ is \textit{regular }(later we'll extend this notion)%
\textit{\ }if the following morphism
\begin{eqnarray}
\widehat{\sigma _{2}} &:&\mathrm{End}\left( \pi \right) \otimes \Sigma
^{1}\rightarrow \mathrm{End}\left( \pi \right) \otimes \Sigma _{1},
\label{connReg} \\
\widehat{\sigma _{2}} &:&\alpha \in \mathrm{End}\left( \pi \right) \otimes
\Sigma ^{1}\longmapsto \left\langle \sigma _{2},\alpha \right\rangle \in
\mathrm{End}\left( \pi \right) \otimes \Sigma _{1},  \notag
\end{eqnarray}%
is an isomorphism.

The next result (cf. \cite{Ly2}) is the direct corollary of the above
theorem.

\begin{thm}
Let $A\in \mathbf{Diff}_{2}\left( \pi ,\pi \right) $ be an operator with
regular symbol. Then, for any connection $\nabla _{c}$ in the cotangent
bundle there is and unique connection $\nabla $ in the bundle $\pi ,$ such
that $\sigma _{1}\left( A,\nabla \right) =0.$
\end{thm}

\section{Symbol classification}

In this section we consider symbols of the second order operators at a point
as tensors of the form $\sigma \in \mathrm{End}\left( E\right) \otimes
S^{2}T,$ where $E$ and $T$ are vector spaces of dimensions $m=\dim E\geq 2\ $%
and$\ \ n=\dim T\geq 2,$ and study their $G=$ $\mathrm{GL}\left( E\right)
\times \mathrm{GL}\left( T\right) $-orbits (cf. \cite{LyB}). All
constructions in this section have pure algebraic nature and can be used over the
field $\mathbb{R}$ as well as the field $\mathbb{C}$ .

We denote by $\sigma _{\theta }\in \mathrm{End}\left( E\right) $ the value
of the symbol on a tensor $\theta \in S^{2}T^{\ast },$ and let
\begin{eqnarray*}
&&g_{\sigma }=(\mathrm{Tr}\otimes \mathrm{id})\sigma \in S^{2}T, \\
&&\left\langle g_{\sigma },\theta \right\rangle =\mathrm{Tr}\left( \sigma
_{\theta }\right) ,
\end{eqnarray*}%
be the quadratic form, associated with the symbol.

Similar to  \cite{LyB} we introduce \textit{Artin-Procesi tensors}%
\begin{equation*}
A_{\sigma ,I}\left( \theta _{1}\otimes \cdots \otimes \theta _{k}\right) =%
\mathrm{Tr}\left( \sigma _{\theta _{i_{1}}}\cdots \sigma _{\theta
_{i_{k}}}\right) ,
\end{equation*}%
where $I=\left( i_{1},...,i_{k}\right) \in S_{k}$ is a permutation of $k$%
-letters, and $\theta _{i}\in S^{2}T^{\ast }.$

Thus, tensors $A_{\sigma ,I}\in \left( \left( S^{2}T^{\ast }\right)
^{\otimes k}\right) ^{\ast }=\left( S^{2}T\right) ^{\otimes k}$ are $G$%
-invariants and $A_{\sigma ,I}=A_{\sigma ,J}$, when permutation $J\in S_{k}$
is obtained from $I$ by a cycle permutation.

Besides Artin-Procesi quadric $g_{\sigma }$, the following tensors $%
h_{2}=A_{\sigma ,\left( 1,2\right) }\in \left( S^{2}T\right) ^{\otimes 2},$ $%
h_{3}=A_{\sigma ,\left( 1,2,3\right) }\in \left( S^{2}T\right) ^{\otimes 3},$%
\begin{eqnarray*}
h_{2}(\theta _{1},\theta _{2}) &=&\mathrm{Tr}\left( \sigma _{\theta
_{1}}\sigma _{\theta _{2}}\right) , \\
h_{3}(\theta _{1},\theta _{2},\theta _{3}) &=&\mathrm{Tr}\left( \sigma
_{\theta _{1}}\sigma _{\theta _{2}}\sigma _{\theta _{3}}\right)
\end{eqnarray*}%
will be extremely important for us.

In what follows we will use \textit{symbols in general position}, or \textit{%
general symbols}, i.e symbols where the following regularity
conditions, together with (4),  hold (cf. \cite{LyB}).

Namely, we require the following:

\begin{enumerate}
\item First of all we require that the quadric $g_{\sigma }\in S^{2}T$ on $%
T^{\ast }$ is nondegenerate. \newline
Denote by $g_{\sigma }^{-1}\in S^{2}T^{\ast }$ the inverse quadric on $T$
and let \ $g_{\sigma }^{\left( 1\right) }=g_{\sigma }^{-1}\rfloor h_{2}\in
S^{2}T,g_{\sigma }^{\left( 2\right) }=g_{\sigma }^{-1}\rfloor g_{\sigma
}^{-1}\rfloor h_{3}\in S^{2}T$ be another quadratic forms on $T^{\ast }$.

\item Denote by $\widehat{g}_{\sigma }^{\left( 1\right) }:T^{\ast
}\rightarrow T^{\ast }$ the linear operator, associated with the pair of
quadratic forms $\left( g_{\sigma },g_{\sigma }^{\left( 1\right) }\right) $
on $T^{\ast }.$ We require that \ $\widehat{g}_{\sigma }^{\left( 1\right) }$
as well as its symmetric power $S=S^{2}\left( \widehat{g}_{\sigma }^{\left(
1\right) }\right) :S^{2}T^{\ast }\rightarrow S^{2}T^{\ast }$are are isomorphisms
with distinct eigenvalues. Moreover, we require that quadrics
\begin{equation*}
g_{\sigma }^{\left( 2\right) },S^{\ast }\left( g_{\sigma }^{\left( 2\right)
}\right) ,.....,\left( S^{\ast }\right) ^{N}\left( g_{\sigma }^{\left(
2\right) }\right) ,
\end{equation*}%
where $N=\dim S^{2}T^{\ast }-1$ are linear independent.

\item The eigenvector basis of the operator $\widehat{g}_{\sigma }^{\left(
1\right) }$ in $T^{\ast }$ (possibly over $\mathbb{C}$ ) we denote by $%
e^{\ast }\left( \sigma \right) $ and the dual basis by $e\left( \sigma
\right) .$ We will require, in addition, that the values $g_{\sigma }\left(
e_{i}^{\ast },e_{i}^{\ast }\right) \neq 0,$ and eigenvectors are normed in
such a way that, $g_{\sigma }\left( e_{i}^{\ast },e_{i}^{\ast }\right) =\pm 1
$.

\item The last regularity condition requires that almost all operators $%
\sigma _{\theta }$ has distinct eigenvalues and the commutator map $\sigma
_{\theta _{1}}\times \sigma _{\theta _{2}}\longmapsto \lbrack \sigma
_{\theta _{1}},\sigma _{\theta _{2}}]$ for these operators is nontrivial.
\end{enumerate}

Denote by $\mathbb{F}$ the splitting field of the characteristic polynomial $%
\det \left( \widehat{g}_{\sigma }^{\left( 1\right) }-\lambda \right) $ over
field $\mathbb{F}_{0}$ of the rational functions on the symbol space $%
\mathrm{End}\left( E\right) \otimes S^{2}T.$

Let $E_{\mathbb{F}}=\mathbb{F}\otimes _{\mathbb{R}}E,$ $T_{\mathbb{F}}=%
\mathbb{F}\otimes _{\mathbb{R}}T$ and let $e^{\ast }\left( \sigma \right) $
let be an eigenvector basis in $T_{\mathbb{F}}^{\ast }$ of the operator $%
\widehat{g}_{\sigma }^{\left( 1\right) }$ and $e\left( \sigma \right) $ be
the dual basis in $T_{\mathbb{F}}..$ Denote by $\lambda \left( \sigma
\right) $ the corresponding eigenvalues.

Let now%
\begin{equation*}
\sigma =\sum_{i,j}\sigma _{ij}\otimes e_{i}\otimes e_{j},
\end{equation*}%
be decomposition of the symbol in the eigenvector basis, where operators $%
\sigma _{ij}\in \mathrm{End}\,E_{\mathbb{F}}$ $.$

Define operators%
\begin{equation*}
R\left( \sigma \right) =\sum_{i,j}\lambda _{i}^{l}\lambda _{j}^{l}\
g_{\sigma }^{\left( 2\right) }\left( e_{i}^{\ast },e_{j}^{\ast }\right)
~\sigma _{ij}\in \mathrm{End}\,E_{\mathbb{F}},
\end{equation*}%
for $l=0,1,...$

Remark that these operators are invariants with respect to transpositions
and therefore all $S_{n}-$permutations of eigenvalues and transformations of
the sign change: $e_{i}\rightarrow \pm e_{i}$.

Therefore, $\ R_{l}\left( \sigma \right) \in \mathrm{End}\,E_{\mathbb{F}_{0}},
$ are operators with rational coefficients.

\begin{lem}
Assume that generality conditions $\left( 1-4\right) $ hold. Then operators $%
R_{l}\left( \sigma \right) ,l=0,1,...,\left( n+2\right) \left( n-1\right) /2,
$ define the symbol $\sigma $ uniquely.
\end{lem}

\begin{pf}
Denote by $\widetilde{\lambda }_{l}$ the following quadratic forms $%
\widetilde{\lambda }_{l}=\sum_{i,j}\lambda _{i}^{l}\lambda
_{j}^{l}~g_{\sigma }^{\left( 2\right) }\left( e_{i}^{\ast },e_{j}^{\ast
}\right) ~e_{i}^{\ast }\otimes e_{j}^{\ast }\in S^{2}\left( T_{\mathbb{F}%
}^{\ast }\right) ,$where $l=0,1,....$

Then $\widetilde{\lambda }_{0}=g_{\sigma }^{\left( 2\right) }\ $and $%
\widetilde{\lambda }_{l+1}=S\left( \widetilde{\lambda }_{l}\right) $.
Therefore, forms $\widetilde{\lambda }_{l},$ $l=0,1,....,\left( n+2\right)
\left( n-1\right) /2$ \ give a basis in $S^{2}\left( T_{\mathbb{F}}^{\ast
}\right) $ and, because of $R_{l}\left( \sigma \right) =\left\langle \sigma ,%
\widetilde{\lambda }_{l}\right\rangle ,$with respect to the natural pairing
$\mathrm{End}\,E_{\mathbb{F}}\otimes S^{2}\left( T_{\mathbb{F}}\right) \times
S^{2}\left( T_{\mathbb{F}}^{\ast }\right) \rightarrow \mathrm{End}\,E_{%
\mathbb{F}},$ operators $R_{l}\left( \sigma \right) $ define the symbol.
\end{pf}

This lemma shows that $G$-equivalence of symbols is equivalent to $\mathrm{%
GL}\left( E_{\mathbb{F}_{0}}\right) $- equivalence of the sets of operators $%
\left\{ R_{l}\left( \sigma \right) \right\} $, and therefore, due to the
Artin-Procesi theorem, is completely described by the Artin-Procesi invariants
\begin{equation}
A_{I}\left( \sigma \right) =\mathrm{Tr}\left( R_{i_{1}}\left( \sigma
\right) \cdots R_{i_{s}}\left( \sigma \right) \right) ,  \label{ArtProInv}
\end{equation}%
where $I=\left( i_{1},...,i_{s}\right) ,\ 0\leq i_{j}\leq \left( n+1\right)
\left( n-2\right) /2.$

The Artin-Procesi theorem (\cite{Pr}) together with the Rosenlicht theorem (%
\cite{Ros}) give us the following result.

\begin{thm}
The field of rational $G$-invariants of symbols $\sigma \in \mathrm{End}\,%
E\otimes S^{2}\left( T\right) $ is generated by the invariants (\ref{ArtProInv}%
). This field separates regular orbits.
\end{thm}

Let's now $\sigma \in \mathrm{End}\,E\otimes S^{2}\left( T\right) $ be a
symbol in general position. Then, due to condition $\left( 3\right) $ the
stationary Lie algebra of the orbit $G\sigma $ is one dimensional, and
therefore
\begin{equation*}
\nu =\mathrm{codim}\left( G\sigma \right) =m^{2}\frac{\left( n+2\right)
\left( n-1\right) }{2}-n^{2}+1.
\end{equation*}

\begin{definition}
In addition to the above notion of regularity we will say that a general
symbol is\textit{\ regular} if there are $\nu $ rational $G$-invariants $%
f_{1},...,f_{\nu }$ that define the orbit $G\sigma $, their differentials
are independent at the points of the orbit, and condition (\ref{connReg})
holds.
\end{definition}

\section{Associated connections and invariants of operators}

Let $\Delta \in \mathbf{Diff}_{2}\left( \pi ,\pi \right) $ be a differential
operator with regular symbol $\sigma =\mathrm{smbl}_{2}\left( \Delta
\right) \in \mathrm{End}\left( \pi \right) \otimes \Sigma _{2}.$ Then the
quadric $g_{\sigma }\in \Sigma _{2}$ is nondegenerate and we denote by $%
\nabla _{c}$ the associated Levi-Civita connection. Taking this connection
in Theorem 2, we get a connection $\nabla =\nabla ^{\Delta }$ in the bundle,
\textit{associated} with operator $\Delta .$ If we denote by $Q_{\Delta }$
the quantization, associated with pair $\left( \nabla ^{\Delta },\nabla
_{c}\right) $ we get the following decomposition of the operator:%
\begin{equation}
\Delta =Q_{\Delta }\left( \sigma \right) +\sigma _{0}\left( \Delta \right) ,
\label{decomp}
\end{equation}%
where $\sigma _{0}\left( \Delta \right) \in \mathrm{End}\left( \pi \right)
, $ is the subsymbol.

Moreover, this decomposition is natural in the sense that
\begin{eqnarray*}
Q_{A_{\ast }\left( \Delta \right) }\left( A_{\ast }\left( \sigma \right)
\right) &=&A_{\ast }\left( Q_{\Delta }\left( \sigma \right) \right) , \\
\sigma _{0}\left( A_{\ast }\Delta \right) &=&A_{\ast }\left( \sigma
_{0}\left( \Delta \right) \right) ,
\end{eqnarray*}%
for any automorphism $A$ of the bundle $\pi .$

Using this construction we introduce invariants of operator to be
Artin-Procesi $\mathrm{GL}\left( \pi \right) $ -invariants for the set of
operators $\left\{ \sigma _{0}\left( \Delta \right) ,R_{k}\left( \sigma
\right) \right\} ,$ where $k=0,1,....,\left( n+2\right) \left( n-1\right)
/2. $

We call them the \textit{Artin-Procesi invariants of the \ operators}.

Remark that they generate rational invariants of the $G$-action on the
\textit{extended symbol space} $\mathrm{End}\left( \pi \right) \otimes
\Sigma _{2}\oplus \mathrm{End}\left( \pi \right) .$ Regular $G$-orbits in
this space has codimension
\begin{equation*}
\nu _{0}=\nu +m^{2},
\end{equation*}%
and $\nu _{0}>n$ as well as $\nu >n.$

Therefore, due to the Rosenlicht theorem, there are $\nu _{0}$ algebraically
independent Artin-Procesi invariants that generate the field of rational $G$%
-invariants on the extended symbol space.

\section{Equivalence of regular differential operators}

We say that a differential operator $\Delta \in \mathbf{Diff}_{2}\left( \pi
,\pi \right) $ is \textit{regular} at a point $p\in M$ \ if its symbol $%
\sigma \in \mathrm{End}\left( \pi \right) \otimes \Sigma _{2}\left(
M\right) $ is regular\ at the point and among Artin-Procesi invariants $%
a_{1},...,a_{\nu _{0}}$ \ of the operator, defining the $G$-orbit of the
pair $\left( \sigma \left( \Delta \right) ,\sigma _{0}\left( \Delta \right)
\right) ,$ there are $n=\dim M$ \ invariants, say $a_{1},..,a_{n},$ such
that differentials of functions $a_{1}\left( \Delta \right) ,..,a_{n}\left(
\Delta \right) $ are linear independent at the point $a\in M.$

Here $a_{i}\left( \Delta \right) $ are the values of invariants at the operator $%
\Delta .$

The functions $a_{i}\left( \Delta \right) ,$ $i=1,..,n,$ define local
coordinates in a neighborhood $U\ni p.$ We call them (\cite{Ly2},\cite{LyB},%
\cite{LYk}) \textit{natural coordinates}.

Taking the values of the basic Artin-Procesi invariants $a_{1},...,a_{\nu
_{0}}$ at $\Delta $ and express them in terms of the natural coordinates $%
a_{i}\left( \Delta \right) ,$ $i=1,..,n,$ we get functions $F_{j},$ such
that
\begin{equation}
a_{j}\left( \Delta \right) =F_{j}\left( a_{1}\left( \Delta \right)
,..,a_{n}\left( \Delta \right) \right) ,\ \ n+1\leq j\leq \nu _{0}.
\label{loc}
\end{equation}

Remark that all rational $G$-invariants of the pair $\left( \sigma \left(
\Delta \right) ,\sigma _{0}\left( \Delta \right) \right) $ are rational
functions of the basic invariant and, due to (\ref{loc}), they are functions
of natural coordinates in the neighborhood.

Similar to (\cite{LyB}) define natural local chart $\phi _{U}:U\rightarrow
\mathbf{D}_{U}\subset \mathbb{R}^{n},$ where%
\begin{equation*}
\phi _{U}\left( b\right) =\left( a_{1}\left( \Delta \right) \left( b\right)
,..,a_{n}\left( \Delta \right) \left( b\right) \right) ,
\end{equation*}%
and let $F_{U}:\mathbf{D}_{U}\rightarrow \mathbb{R}^{\nu _{0}-n}$ be the
function given by (\ref{loc}).

We call the triple $\left( \phi _{U},\mathbf{D}_{U},F_{U}\right) $ \textit{%
the model of }$\Delta $ in $U.$

\begin{thm}
Let differential operators $\Delta $,$\Delta ^{\prime }\in \mathbf{Diff}%
_{2}\left( \pi ,\pi \right) $ have the same model in the open set $U.$
\newline
Then there is an automorphism $A_{U}\in \mathrm{GL}\left( \pi _{U}\right) $
of the restriction of the bundle $\pi $ on\ domain $U$ such that $A_{U\ast }\left(
\Delta \right) =\Delta ^{\prime }.$
\end{thm}

\begin{pf}
Condition that operators have the same model means that the pairs $\left(
\sigma \left( \Delta \right) ,\sigma _{0}\left( \Delta \right) \right) $ and
$\left( \sigma \left( \Delta ^{\prime }\right) ,\sigma _{0}\left( \Delta
^{\prime }\right) \right) $ belong to the same $G$-orbit at any point $a\in
U.$ Therefore, there exists the above automorphism $A_{U}.$ This means that
the quadratic forms $g_{\sigma }$ coincide in $U$ and we have $\nabla
^{\prime }=A_{U\ast }(\nabla )$ for the associated connections. Therefore, $%
A_{U\ast }\left( \Delta \right) =\Delta ^{\prime }$ in the neighborhood.
\end{pf}

\begin{corollary}
Let differential operators $\Delta $,$\Delta ^{\prime }\in \mathbf{Diff}%
_{2}\left( \pi ,\pi \right) $ has models $\left( \phi _{U},\mathbf{D}%
_{U},F_{U}\right) $ and $\left( \phi _{U^{\prime }},\mathbf{D}_{U^{\prime
}},F_{U^{\prime }}\right) $ in open sets $U\subset M$ and $U^{\prime
}\subset M$ \ defined by the same basic invariants $a_{1},..,a_{n}$ and
functions $F_{j}.$ \newline
Let $\widetilde{U}\subset U$ and $\widetilde{U^{\prime }}\subset U^{\prime }$
be open and simply connected domains such that $\phi _{U}\left( \widetilde{U}%
\right) =\phi _{U^{\prime }}\left( \widetilde{U^{\prime }}\right) \subset
\mathbf{D}_{U}\cap \mathbf{D}_{U^{\prime }}.$ \newline
Then there is an automorphism $A_{\widetilde{U},\widetilde{U^{\prime }}}:\pi
_{\widetilde{U}}\rightarrow \pi _{\widetilde{U^{\prime }}}$ \ such that $%
\left( A_{\widetilde{U},\widetilde{U^{\prime }}}\right) _{\ast }\Delta
=\Delta ^{\prime }.$
\end{corollary}

Following (\cite{LyB}) we say that an atlas $\left\{ \left( \phi _{U^{\alpha
}},\mathbf{D}_{U^{\alpha }}\right) \right\} $ given by models $\left\{ \phi
_{U^{\alpha }},\mathbf{D}_{U^{\alpha }},F_{U^{\alpha }}\right\} $ is \textit{%
natural} \ if the\ sets of basic invariants $\left( a_{1}^{\alpha
},...,a_{n}^{\alpha }\right) $ are different for different $\alpha .$

Let the operators $\Delta $ and $\Delta ^{\prime }$ have the same natural atlas.
Then isomorphisms $A_{U_{\alpha }}$ as well as $A_{U_{\alpha }\cap U_{\beta
}}$ are defined up to multiplication by nonvanishing functions $f_{\alpha }$
and $f_{\alpha \beta }.$ Therefore, the existence of the global isomorphism
is equivalent to triviality of the cohomology class in $H^{1}\left( M,%
\mathbb{Z}_{2}\right) ,$ defining by the cocycle $\left\{ f_{\alpha \beta
}\right\} .$ Summarizing we get the following.

\begin{thm}
Two regular linear differential operators $\Delta ,\Delta ^{\prime }\in
\mathbf{Diff}_{2}\left( \pi ,\pi \right) $ on a manifold $M,$ with $%
H^{1}\left( M,\mathbb{Z}_{2}\right) =0,$ are $\mathrm{Aut}\left( \pi
\right) $-equivalent if and only if a natural atlas for operator $\Delta $
is the natural atlas for $\Delta ^{\prime },$ i.e. if they have the same
models.
\end{thm}

\section*{Acknowledgements}
The author was partially supported by the Russian Foundation for Basic
Research (project 18-29-10013).


\end{document}